\documentclass[11pt]{amsart}
\usepackage{geometry,amsthm,graphicx,float,enumerate}
\usepackage{amsmath}
\usepackage[makeroom]{cancel}
\usepackage{scalefnt}
\usepackage{mathrsfs}
\usepackage{color}
\usepackage{subfiles}
\usepackage{enumitem}

\usepackage{charter}
\usepackage{euler}
\usepackage{amssymb}
\restylefloat{table}
\restylefloat{figure}

\swapnumbers
\newtheorem{thm}{Theorem}[section]

\newtheorem{cor}[thm]{Corollary}

\newtheorem{lemma}[thm]{Lemma}
\newtheorem{prop}[thm]{Proposition}
\theoremstyle{definition}
\newtheorem{defn}[thm]{Definition}
\theoremstyle{remark}
\newtheorem{rem}[thm]{Remark}

\newtheorem{prob}[thm]{Problem}

\newtheorem{solution}[thm]{Solution}

\newtheorem{notn}[thm]{Notation}

\newcommand{\BB}[1]{\mathbb{#1}}
\newcommand{\CAL}[1]{\mathcal{#1}}

\newcommand{\INDEX}[1]{\big[{[}#1{]}\big]}

\newcommand{\DF}{d_{\BB{F}^m}}
\newcommand{\KF}{k_{\BB{F}^m}}

\newcommand{\TR}{\text{tr}}
\newcommand{\SPAN}{\text{span}}

\newcommand{\IM}{\text{im}}

\newcommand{\BLOCKSET}{\BB{S}}
\newcommand{\GRMAN}[3]{\mathcal G(#1, \mathbb  {#2}^{#3} )}
\newcommand{\GRASS}{\GRMAN{l}{\mathbb F}{m}}
\newcommand{\GRMIN}{\mu_{n,l,\mathbb F^m}}

\newcommand{\CBRACK}[1]{ \left\{ #1 \right\} } 
\newcommand{\PARENTH}[1]{ \left( #1 \right) } 

\makeatletter
\providecommand*{\cupdot}{%
  \mathbin{%
    \mathpalette\@cupdot{}%
  }%
}

\newcommand*{\@cupdot}[2]{%
  \ooalign{%
    $\m@th#1\cup$\cr
    \sbox0{$#1\cup$}%
    \dimen@=\ht0 %
    \sbox0{$#1\cdot$}%
    \advance\dimen@ by -\ht0 %
    \dimen@=.5\dimen@
    \hidewidth\raise\dimen@\hbox{$\m@th#1\cdot$}\hidewidth
  }%
}

\begin{document}

\title[Optimally spread subspace packings]{Constructions and properties of optimally spread subspace packings via symmetric and affine block designs and mutually unbiased bases
}

\author{Peter G. Casazza}
\author{John I. Haas IV}
\author{Joshua Stueck}
\author{Tin T. Tran}



\thanks{All authors of this paper were supported by NSF DMS 1609760, NSF DMS 1725455 and ARO W911NF-16-1-0008.}

\address{219 Mathematical Sciences Building, Department of Mathematics, University of Missouri, Columbia, MO 65211}

\begin{abstract}
We continue the study of optimal chordal packings, with emphasis on packing subspaces of dimension greater than one.  Following a principle outlined in ~\cite{2016arXiv160704546B}, where the authors use maximal affine block designs and maximal sets of mutually unbiased bases to construct Grassmannian $2$-designs, we show that their method extends to other types of block designs, leading to a plethora of optimal subspace packings characterized by the orthoplex bound.  More generally, we show that any optimal chordal packing is necessarily a fusion frame and that its spatial complement is also optimal.
\end{abstract}
$\,\,$
\\

\maketitle


\section{Introduction}

At a recent AMS meeting on Advances in Packings~\cite{ams_pack_2018} at Ohio State University, several participants agreed that a seemingly inherent  number-theoretic principle underpins the existence of most, if not all, Symmetric Informationally Complete Positive Operater Valued Measures, or simply {\it SICs} -- a special type of chordal packing in complex projective space, $\mathbb C \mathbb P^{m-1}$, consisting of $m^2$ lines with constant pairwise chordal distance between its elements -- objects which are of particular interest to quantum information theorists~\cite{Zauner1999,Appleby2009, 2017arXiv170303993S}. Accordingly, some speakers advocated that -- by developing a tabulation of analytic and numerical constructions of SICs for different parameters (ie, for each dimension, $m$) -- the community might eventually extrapolate an existential, or even constructive, proof of Zauner's conjecture, which states that SICs exist in every dimension.

This work is motivated by a similar but generalized principal. By developing a parametric tabulation for all optimally spread packings with respect to the chordal distance in arbitrary Grassmannian manifolds, we hope to contribute to a fundamental theorem for the construction and structure of all optimally spread packings. 

With regard to developing a universal theory of said structure, in Section~\ref{sec:prop}, we prove that all optimally spread packings are fusion frames.  In addition, we show that their spatial complements are also optimally spread.   

Toward the development of the aforementioned tabulation of solutions, we demonstrate elementary constructions of optimally spread subspace packings for parameter sets characterized by the simplex bound; less trivially, we construct optimal packings characterized by the orthoplex bound (see Section~\ref{sec:recipe} for information about both bounds) by exploiting a technique presented in~\cite{2016arXiv160704546B} which led to a special class of both real and complex solutions, so-called {\it maximal orthoplectic fusion frames}, which had previously been constructed for the real case in \cite{ShorSloane1998} by other means.  

In this vein, we re-present an abbreviated proof of the method of construction from~\cite{2016arXiv160704546B}, emphasizing its dependence on the existence of certain families of mutually unbiased bases and certain block designs, along with a few other insights.  Finally, we exploit the existence of numerous other classes of block designs which -- omitted in~\cite{2016arXiv160704546B} -- which, to our surprise also  satisfy the restrictive conditions of the construction, thereby by producing more optimally spread fusion frames characterized by the orthoplex bound.

\section{Preliminaries}\label{sec:prelim}


Throughout, we assume that $l,m,$ and $n$ are positive integers satisfying $l \leq m \leq n$ and that  $\mathbb F= \mathbb R$ or $\mathbb F= \mathbb C$.  Recall that the {\bf Grassmannian manifold}, denoted $\GRMAN{l}{F}{m}$, is the space of all $l$-dimensional subspaces in the Hilbert space, $\mathbb F^m$.
We refer to a finite sequence of subspaces, $\CBRACK{ \mathbb W_j}_{j=1}^n \subset \GRMAN{l}{F}{m}$, as an {\bf $(n,l)$-packing for $\mathbb F^m$}, or simply a {\bf packing} when the context is clear.

As a minor abuse of notation, we interchangeably refer to a given packing $\mathcal P$ either by its sequence of subspaces, $\mathcal P = \CBRACK{ \mathbb W_j}_{j=1}^n$, or its corresponding sequence of unique orthogonal projections, $\mathcal P = \CBRACK{P_j}_{j=1}^n$:  
 \begin{center}ie, $\IM(P_j)=\mathbb W_j$, $\TR(P_j)=l$, and $P_j=P_j^2=P^*$ for every $j \in \CBRACK{1, \dots, n}$.  
 \end{center} 
	To be precise, in Section~\ref{sec:prop}, it is convenient to alternate between these equivalent interpretations, while in the remaining sections we primarily identify packings as sequences of orthogonal projections. 

As considered in~\cite{ConwayHardinSloane1996, MR2915706}, there are numerous notions of distance that one may define between two subspaces.  This work concerns the study and construction of packings which are ``optimally spread'' with respect to  
the {\it chordal distance}.

\begin{defn}
	Given two  $l$-dimensional subspaces of $\BB{F}^m$ with corresponding orthogonal projections $P$ and $P'$, the {\bf chordal distance} between them is
	$$
	d_c(P,P') := \frac{1}{\sqrt 2} \| P - P'\| = \PARENTH{l - \TR( P P' ) }^{1/2} \, . 
	$$  
\end{defn}

With respect to this objective function, our goal is to study and construct packings which maximize the minimal pairwise distance over the space of all $n$-packings in $\GRMAN{l}{F}{m}$.  To facilitate this, 
we define and denote the {\bf (chordal) coherence} of an $(n,l)$-packing, $\mathcal P =\{P_j\}_{j=1}^n$, for $\mathbb F^m$ as
$$
\mu(\mathcal P) := \max_{\substack{1\le j,j'\le n\\ j \neq j'}} \TR(P_j P_{j'})
$$ 
and the {\bf packing constant} as
$$
\GRMIN := \inf_{\mathcal P \subset \GRMAN{l}{F}{m}} \mu(\mathcal P).
$$ 
Formally, we say an $(n,l)$-packing, $\mathcal P$, for $\mathbb F^m$ is {\bf optimally spread} if
$$
\mu(\mathcal P) = \GRMIN.
$$
An elementary topological argument ensures that optimally spread packings exist for all parameters satisfying $l \le m \le n$.
\begin{prop}\label{prop:optimalexists}
	An optimally spread $(n,l)$-packing for $\mathbb F^m$ exists.
\end{prop}
\begin{proof}
	The Grassmannian manifold, $\GRMAN{l}{F}{m}$, is well known to be compact; whence, the space of all $(n,l)$-packings for $\mathbb F^m$ -- which can be identified with the Cartesian product $\Pi_{j=1}^{n} \GRMAN{l}{F}{m}$ -- is also compact.  Because the coherence function is continuous, the claim follows by the extreme value theorem.
\end{proof}

These definitions and the preceding proposition imply that the following problem well-posed. Formally, the {\it (constant-rank) chordal packing problem} is stated as follows.

\begin{prob}\label{prob:chord}[{\bf The (constant-rank) chordal packing problem}]
For each quintuplet of parameters, $(n,l,m, \mathbb F)$, determine the corresponding packing constant, $\GRMIN$, and, if possible, construct and characterize an optimally spread $n$-packing, $\mathcal P \subset \GRMAN{l}{F}{m}$.
\end{prob}

From an applications point of view, we are especially interested in $(n, l)$-packings endowed with additional spectral constraints, which engender numerous signal processing possibilities~\cite{MR2964005, Bachoc2004, 5946867, MR2860105, bod_cas_edi_bal_2008, MR1984549, MR1829801, MR3367817, MR2399108, MasseyRuizStojanoff2010, BachocEhler2013, RoyScott2007, chebira2008frames}; in particular, we are interested in fusion frames.

\begin{defn}\label{def:frames}
	An $(n,l)$-packing for $\mathbb F^m$, $\CAL{P}=\{P_j\}_{j=1}^n$, is an {\bf $(n,l)$-fusion frame} for $\mathbb F^m$ if the projections' sum, $\sum_{j=1}^n P_j$, is positive definite, and whenever $l=1$, then $\mathcal P$ is also called a {\bf frame}. Any $(n,l)$ fusion frame   is {\bf tight} if this sum satisfies 
	 $$\sum_{j=1}^n P_j = \frac{nl}{m} I_m,$$ where -- henceforth -- $I_m$ denotes the $m \times m$ identity matrix. An $(n,l)$-fusion frame for $\mathbb F^m$ is {\bf Grassmannian} if it is optimally spread.
\end{defn}

\begin{rem} $\,$ ~\\
	\begin{enumerate}
		\item In the case $l=1$, where a packing is a frame, the projectors are  usually - up to a choice of unimodular phasing - identified with unit vectors.
		\item 
		      The term ``Grassmannian fusion frame'' is an acknowledgement of~\cite{MR1984549}, where the authors were interested specifically in frames (ie, the case $l=1$), or equivalently packings in projective space where there is no ambiguity in the use of the term ``Grassmannian"; however, because we are now studying packings in arbitrary Grassmannian manifolds, we prefer the term ``optimally spread'' over ``Grassmannian'' in order to avoid confusion when describing solutions to Problem~\ref{prob:chord}.
	\end{enumerate}
\end{rem}

While it is false that all $n$-packing are fusion frames, fortunately for applications - as we show next - all optimally spread $(n,l)$-packings are fusion frames with the corresponding parameters.

\section{Some properties of optimal packings}\label{sec:prop}
In this section, we prove two important facts about optimal packings.  The first concerns their spanning properties.

\subsection{All optimal packings are fusion frames} 

It is known~\cite{2017arXiv170701858F} that if an $(n,1)$-packing for $\mathbb F^m$ is an optimally spread packing, then it is a Grassmannian frame, but we are unaware of an analogous statement for the general case of optimal $(n,l)$-packings, $1 \le l \le m$.  In the following, we show that all optimal packings are fusion frames.
To prove this, we begin with two lemmas.

\begin{lemma}\label{lem_span1}
Suppose $\GRMIN, \delta$ and $\epsilon$ are positive real numbers satisfying $\GRMIN>\epsilon$
and
\begin{equation}\label{eq:lemspan1}
 2\sqrt{l\delta(\GRMIN-\epsilon)} + l\delta < \frac \epsilon 2.
\end{equation}
If
 $\CBRACK{x_j}_{j=1}^l, \CBRACK{y_j}_{j=1}^l$ and $\CBRACK{z_j}_{j=1}^l$  are sequences of unit vectors in $\mathbb F^m$ satisfying
\begin{equation}\label{eq:lemspan2}
\sum_{j=1}^l \sum_{k=1}^l \left| \langle x_j, y_k\rangle \right|^2 < \GRMIN - \epsilon
\end{equation}
and 
\begin{equation}\label{eq:lemspan3}
\sum_{j=1}^l \|z_j - y_j\|^2 < \delta,
\end{equation} then
\begin{equation*}
\sum_{j=1}^l \sum_{k=1}^l \left| \langle x_j, z_k\rangle \right|^2 < \GRMIN - \frac \epsilon 2.
\end{equation*}
\end{lemma}

\begin{proof}
	Rewriting each $z_j$ as $z_j=y_j +(z_j-y_j)$ yields
	$$
	\sqrt{\sum_{j=1}^l \sum_{k=1}^l \left| \langle x_j, z_k\rangle \right|^2} \leq
	\sqrt{\sum_{j=1}^l \sum_{k=1}^l \left| \langle x_j, y_k\rangle \right|^2}
	+
	\sqrt{\sum_{j=1}^l \sum_{k=1}^l \left| \langle x_j, z_k - y_k\rangle \right|^2}.
	$$
	Next, the Cauchy-Schwarz inequality and (\ref{eq:lemspan2}) imply
	$$
	\sqrt{\sum_{j=1}^l \sum_{k=1}^l \left| \langle x_j, z_k\rangle \right|^2}
	\leq
	\sqrt{\GRMIN - \epsilon} + \sqrt{\sum_{j=1}^l \sum_{k=1}^l \|x_j\|^2 \|z_k-y_k\|^2} 
	<
	\sqrt{\GRMIN - \epsilon} + \sqrt{l \delta},
	$$
	where the last inequality in the previous line follows from Equation~(\ref{eq:lemspan3}) and the unit norm property of the vectors.  The claim follows by squaring and applying Equation~(\ref{eq:lemspan1}) to obtain
	$$
	\PARENTH{\sqrt{\GRMIN - \epsilon} + \sqrt{l \delta}}^2 = \GRMIN - \epsilon + 2\sqrt{l\delta(\GRMIN-\epsilon)} + l\delta < \GRMIN - \epsilon + \frac \epsilon 2 = \GRMIN - \frac \epsilon 2.
	$$
\end{proof}

To facilitate the next lemma and the theorem that follows, we say an element $P_k$ (or its image, $\mathbb W_k$) of an optimally spread $(n,l)$-packing, $\mathcal P=\CBRACK{P_j}_{j=1}^n$, for $\mathbb F^m$ {\bf achieves the packing constant} if there exists $k' \in \{1,2,...,n\}$ with $k \neq k'$ such that $\TR(P_k P_{k'})=\GRMIN$, and we call each element $P_{k'}$ satisfying this condition a {\bf packing neighbor} of $P_k$.

\begin{lemma}\label{lem_span2}
	Suppose $\mathcal P=\CBRACK{P_j}_{j=1}^n$ is an optimally spread $(n,l)$-packing in $\mathbb F^m$ with corresponding subspaces $\CBRACK{\mathbb W_j}_{j=1}^n$ and, furthermore, suppose that $\GRMIN>0$. If $\mathbb W_k$ is an element of $\mathcal P$ that achieves the packing constant, then it contains a unit vector  which is not orthogonal to any of its packing neighbors.
\end{lemma}

\begin{proof}
	Write 
	$\mathcal I:=\CBRACK{1\leq i \leq n: \TR(P_i P_k) = \GRMIN}$,
	and for each $i \in \mathcal I$, let $\mathbb V_i$ denote the maximal subspace of $\mathbb W_k$ which is orthogonal to $\mathbb W_i$, ie $\mathbb V_i := \mathbb W_k \cap \mathbb W_i^\perp$.  The assumption $\GRMIN>0$ implies that every $\mathbb V_i$ is a proper subspace of $\mathbb W_k$, and since a linear space cannot be written as a finite union of proper subspaces, it follows that $\mathbb W_k \backslash \cup_{i \in \mathcal I} \mathbb V_i$ is nonempty, so the claim follows.
\end{proof}

\begin{thm}\label{thm:span}
	If $\mathcal P=\CBRACK{P_j}_{j=1}^n$ is an optimally spread $(n,l)$-packing in $\mathbb F^m$ with corresponding subspaces $\CBRACK{\mathbb W_j}_{j=1}^n$ and 
	$$\mathcal K := \CBRACK{1 \leq k \leq n : \mathbb W_k \text{ achieves the packing constant}},$$ then
	$\SPAN\CBRACK{\mathbb W_k : k \in \mathcal K} = \mathbb F^m,$ meaning $\CBRACK{\mathbb W_k}_{k \in \mathcal K}$ is an $(n',l)$-fusion frame for $\mathbb F^m$, where $n':=\left| \mathcal K \right|$.
	
Consequently, $\mathcal P$ is an $(n,l)$-fusion frame for $\mathbb F^m$ and, in particular, a Grassmannian fusion frame.
\end{thm}
    
\begin{proof}
If $\GRMIN=0$, then $\mathcal K = \CBRACK{1,2,\dots,n}$, so $\mathbb W_j$ is in the orthogonal complement of $\mathbb W_{j'}$ for every $j \neq j'$, and since $nl\geq m$, it follows that $\SPAN\CBRACK{\mathbb W_j}_{j=1}^n = \mathbb F^m$.	For the case $\GRMIN>0$, we proceed by way of contradiction, iteratively replacing elements of $\mathcal P$ that achieve the packing constant in such a way that we eventually obtain a new $(n,l)$-packing in $\mathbb F^m$ with coherence strictly less than $\GRMIN$, which cannot exist.
With the contradictory approach in mind, fix a unit vector $z \in \mathbb F^m$ so that $z$ is in the orthogonal complement of $\SPAN\CBRACK{\mathbb W_j}_{j \in \mathcal K}$.
Next, we describe the replacement procedure.

The first step is to choose some $k \in \mathcal K$, write
$$
\mathcal I_k:= \CBRACK{1\leq i \leq n : \TR(P_i P_k) = \GRMIN}
\text{ and }
\mathcal I_k^c:= \CBRACK{1\leq i \leq n : \TR(P_i P_k) < \GRMIN},
$$
and for every $i \neq k$, fix an orthonormal basis $\{x_{i,j}\}_{j=1}^l$ for $\mathbb W_i$. By Lemma~\ref{lem_span2}, there exists a unit vector $x_{k,1} \in \mathbb W_k$ which is nonorthogonal to all of $\mathbb W_k$'s packing neighbors, so apply the Gram-Schmidt algorithm to extend $x_{k,1}$ to an orthonormal basis, $\{x_{k,j}\}_{j=1}^l$, for $\mathbb W_k$.  

As defined, there exists $\epsilon>0$ so that
$$
 \sum_{j=1}^l \sum_{j'=1}^l 
 \left| 
 \langle  x_{i,j}, x_{k,j'} \rangle 
 \right|^2 = \TR(P_i P_k)  < \GRMIN - \epsilon, \text{ for all } i \in \mathcal I_k^c.
$$
Choose $\delta \in \mathbb R$ so that $0 < \delta <1/2$
and (\ref{eq:lemspan1}) from Lemma~\ref{lem_span1} is satisfied, define
$$
y_{k,1}:=\sqrt{1-\delta^2} \, x_{k,1} + \delta z \text{ and } y_{k,j}:=x_{k,j} \text { for } 2 \leq j \leq l,
$$
and define $\mathbb V_k := \SPAN\CBRACK{y_{k,j}}_{j=1}^l$ with corresponding orthogonal projection, $Q_k$.  Because $\langle x_{k,1}, z \rangle =0$, it follows by elementary computation that $\CBRACK{y_{k,j}}_{j=1}^l$ is an orthonormal basis for $\mathbb V_k$.

Noting that $1-\delta^2 < \sqrt{1-\delta^2}$ implies $2\PARENTH{1 - \sqrt{1-\delta^2}} < 2 \delta^2$, we estimate
\begin{equation*}
\sum_{j=1}^l \|x_{k,j} - y_{k,j} \|^2 
= 
\|x_{k,1} - y_{k,1} \|^2 
= 2 \PARENTH{1 - \sqrt{1-\delta^2}} < 2\delta^2 < \delta.
\end{equation*}
Therefore, Lemma~\ref{lem_span1}
  implies
$$
\TR(P_i Q_k) = \sum_{j=1}^l \sum_{j'=1}^l \left| \langle x_{i,j}, y_{k, j'}  \rangle  \right|^2 < \GRMIN - \frac \epsilon 2 \text{ for all } i \in \mathcal I_k^c.
$$
Furthermore, the nonorthogonality of $x_{k,1}$ with $\mathbb W_i$ for every $i \in \mathcal I_k$ implies
\begin{equation*}
0<
\sum_{j=1}^l \left|\langle y_{k,1}, x_{i,j}  \rangle\right|^2 = (1 - \delta^2) \sum_{j=1}^l \left|\langle x_{k,1}, x_{i,j}  \rangle\right|^2 < \sum_{j=1}^l \left|\langle x_{k,1}, x_{i,j}  \rangle\right|^2  \text{ for all } i \in \mathcal I_k,
\end{equation*}
so it follows that
$$
\TR(P_i Q_k) 
=
 \sum_{j=1}^l \sum_{j'=1}^l \left|\langle   y_{k,j}, x_{i,j'} \rangle \right|^2
 <
  \sum_{j=1}^l \sum_{j'=1}^l \left|\langle   x_{k,j}, x_{i,j'} \rangle \right|^2
  =
  \TR(P_k P_i) = \GRMIN
  \text{ for all } i \in \mathcal I_k.
$$
Thus, replacing $\mathbb W_k$ with $\mathbb V_k$ produces a new $(n,l)$-packing for $\mathbb F^m$,  where the replaced element no long achieves the packing constant.  

Now, we iterate this replacement procedure. After at most $n'$ repetitions of this process, we  obtain a final $(n,l)$-packing, 
with coherence strictly less than $\GRMIN$, the desired contradiction, so the claims follow.
\end{proof}

Next, we consider the {\it spatial complements} of Grassmannian frames. 

\subsection{Spatial complements of Grassmannian fusion frames}
Given an $(n,l)$-packing, $\mathcal P = \CBRACK{P_j}_{j=1}^n$, for $\mathbb F^m$, its {\bf spatial complement} is the $(n, m-l)$-packing,  $\mathcal P^\perp := \CBRACK{I_m - P_j}_{j=1}^n$,  where $I_m$ denotes the $m \times m$ identity matrix.  As one might expect, the spatial complement of a Grassmannian fusion frame is also optimally spread, which we show in the following theorem.

\begin{thm}
	If $\mathcal P=\CBRACK{P_j}_{j=1}^n$ is a Grassmannian $(n,l)$-fusion frame for $\mathbb F^m$, then its spatial complement, $\mathcal P^\perp$, is a Grassmannian $(n, m-l)$-fusion frame for $\mathbb F^m$ and $\mu_{n,m-l,\mathbb F^m} = m - 2l + \GRMIN.$
\end{thm}
\begin{proof}
	Given any $(n,l)$-packing for $\mathbb F^m$, say $\mathcal Q:=\{Q_j\}_{j=1}^n$, we compute the coherence of its spatial complement, $\mathcal Q^\perp$, 
	$$
	\mu\PARENTH{\mathcal Q^\perp} 
	=
	\max_{j \neq j'} \TR\PARENTH{(I_m - Q_j)(I_m - Q_{j'})} 
	=
	m - 2l + \max_{j \neq j'} \TR\PARENTH{Q_j Q_{j'}}
	=
	m - 2l + \mu\PARENTH{\mathcal Q},
	$$
	thereby showing that the coherence of $\mathcal Q^\perp$ depends only on the coherence of $\mathcal Q$.
	Thus, $\mu\PARENTH{\mathcal P^\perp} = m - 2l + \GRMIN$ is minimal over the space of all $(n,m-l)$-packings for $\mathbb F^m$, so $\mathcal P^\perp$ is optimally spread and therefore a Grassmannian $(n,m-l)$-fusion frame for $\mathbb F^m$ by Theorem~\ref{thm:span}. 
\end{proof}

\section{Recalling a recipe}\label{sec:recipe}

With an equivalence between optimally spread $n$-packings and Grassmannian fusion frames of corresponding parameters established, this section is dedicated to the construction of infinite families of tight, optimally spread fusion frames.  
Our constructions rely on three things:
\begin{itemize}
	 \item[(I)]
the most successfully obtained lower chordal coherence bounds~\cite{WoottersFields1989, Fickus:2015aa, MR3557826, BandyopadhyayBoykinRoychowdhuryVatan2002} -- the simplex and orthoplex bounds -- often attributed to Rankin~\cite{Rankin1955} and Conway, Hardin, and Sloane~\cite{ConwayHardinSloane1996}, and Welch~\cite{Welch1974},  which motivate ours study,  
     \item[(II)] two major ingredients, the existence of maximal sets of mutually unbiased bases in certain settings~\cite{WoottersFields1989, BoykinSitharamTarafiWocjan} and the existence of certain $1$-block designs~\cite{MR2246267}, and 
     \item[(III)] the technique, presented in~\cite{2016arXiv160704546B}, which combines the two ingredients to reveal numerous infinite families of optimally spread packings.
\end{itemize}
  
\subsection{Two coherence bounds}
     The {\it simplex} and {\it orthoplex} bounds admit several manifestations, for example as seen in~\cite{Rankin1955, Welch1974, ConwayHardinSloane1996, Levenshtein1992, BargNogin2002}. Of relevance to this work is their occurrence (i) in the {\it coding problem}, where one maximizes the minimal distance between a prescribed number of points on a real sphere of fixed dimension~(see~\cite{Rankin1955, Ericson2001, BargNogin2002, 2006math.6734B}), (ii)  its connection with a {\it constrained version of the coding problem} induced by the chordal packing problem via the so-called {\bf $l$-tracelesss map}, henceforth denoted $T_l$, which seemingly laid the foundation for the mainstream chordal packing arena, due to Conway, Hardin, and Sloane in 1996~\cite{ConwayHardinSloane1996}, and (iii) the manner with which the constrained coding problem is equivalent to the chordal packing problem. 
     
     Given $1 \le l \le m$, the {\bf$l$-traceless map}  is
  \begin{equation}\label{eq:dimtr0}
     T_l:\GRMAN{l}{\mathbb F}{m} \rightarrow \mathbb R^{\DF} : P \mapsto \mathcal  V\left(P-\frac l m I_m\right) \in \mathcal V(\mathbb H),
   \end{equation}
      where $\mathbb H $ denotes the corresponding $m\times m$ matrix subspace of trace zero elements in $\mathbb F^{m \times m}$, $\mathcal V:\mathbb H \rightarrow \mathbb R^{\DF}$ is a fixed vectorization isomorphism that maps the traceless symmetric/hermitian matrices of $\mathbb H$ to vectors in $\mathbb R^{\DF}$, and  where the vanishing trace constraint implicates the two isomorphic spaces' real dimension,
      \begin{equation}\label{eq:DF}
  \DF := \left\{ \begin{array}{cc} \frac{(m+2)(m-1)}{2}, & \BB{F}= \BB{R} \\  m^2 -1, & \BB{F} =\BB{C} \end{array} \right.  .
  \end{equation}
  
      \begin{notn}
      Henceforth, as computed above, we denote by $\DF$ the  real dimension of the so-called ``traceless' space'', $\mathbb H$, or equivalently, the dimension of real vectorized space, $\mathcal V(\mathbb H)$ into which an element of $\GRASS$ embeds via the $l$-traceless map.
      \end{notn}

     It is elementary to verify that $T_l$ is a  (scaled) isometry~\cite{ConwayHardinSloane1996, bgb15}; in particular,  for elements $P, P' \in \GRASS$, the $l$-traceless map enforces the {\bf traceless identity}:   
   \begin{equation}\label{eq:tr0_id}
      \TR\PARENTH{PP'} = \frac {l^2}{m} + \frac{l(m-l)}{m}  \left\langle v_P^{(l)}, v_{P'}^{(l)} \right\rangle,
   \end{equation}
     where the embedded, rescaled-to-unit vectors are $$v_P^{(l)}:=\frac{T_l(P)}{\|T_l(P)\|}, \quad v_{P'}^{(l)}:=\frac{T_l(P')}{\|T_l(P')\|} \in \mathcal S^{\DF-1} \subset \mathbb R^{\DF},$$ and the vectorization, $\mathcal V$, converts the trace inner product between points in the matrix subspace $\mathbb H$ into the standard inner product between points in $\mathbb R^{\DF}$.
     \begin{notn}
     	As in the last line and  Equation~\ref{eq:tr0_id}, if no other indexing scheme is established -- which will occur occassionally in this work -- projections embedded via $T_l$ are sub-indexed by their underlying projection, and - although unnessary - it is convenient to record the underlying projections' ranks in the superscripts.
     \end{notn} 
 
     Of significance here and in future work (where we are studying optimally spread mixed rank packings), placing the appropriate ``$\min \max$'' statements in front of the left and right hand sides of Equation~\ref{eq:tr0_id} converts Problem~\ref{prob:chord}, the chordal optimization problem, into a {\bf restricted coding problem}.    
     
     \begin{defn}
     Given any $d,n \in \mathbb N$, an {\bf $n$-code}, $\mathcal C$, is a sequence of $n$ unit vectors on the real unit sphere in $\mathbb R^d$, ie $\mathcal C:=\{v_i\}_{i=1}^n \subset \mathcal S^{d-1} \subset \mathbb R^d$.
     \end{defn}	
     
     \begin{prob}\label{prob:restr}[{\bf The restricted coding problem}]
     	Let $\Omega:=T_l(\GRASS) \subset \mathcal S^{\DF-1} \subset \mathbb R^{\DF}$, the normalized image of $\GRASS$ under the $l$-traceless map.
     	Determine 
     	$$
     	\sigma_{n,l, \mathbb F^m} : = \inf_{\mathcal C:=\CBRACK{v_i}_{i=1}^n\subset \Omega } \max_{i \neq j} \langle v_i, v_j \rangle
        $$
     \end{prob} 
  
     \begin{rem}
     	As with the chordal packing problem, we reserve the symbol $\sigma_{n,l, \mathbb F^m}$ to refer to the solution to this problem.
     \end{rem}

     Because the $l$-traceless map is continuous, an elementary topological proof similar to that of Proposition~\ref{prop:optimalexists} verifies the existence of solutions for all parameters in the aforementioned problem.  Perhaps the more obvious proof is that the chordal packing problem was already shown to be  well-defined, so its  equivalence to Problem~\ref{prob:restr} via the traceless relationship assures well-definedness.  Indeed, any solution to Problem~\ref{prob:restr}  resolves to a corresponding solution for Problem~\ref{prob:chord}.
     
     \begin{solution}[General solution]
     	 For every quadruple of parameters, $(n,l,m,\mathbb F)$,
     	  \begin{equation}\label{eq:gensol} \GRMIN = \frac {l^2}{m} + \frac{l(m-l)}{m}  	\sigma_{n,l, \mathbb F^m}.\end{equation}
     \end{solution}
 
     Unfortunately,  few  -- in fact, finitely many~\cite{lev:2017}, except in $\mathbb R^2$~\cite{BK06}) -- solutions for this plethora of problems are known which are characterized neither by the simplex nor orthoplex bounds, which -- incidentally -- represent a complete set of tight bounds for the {\bf unrestricted coding problem}, under suitable conditions.

       \begin{prob}\label{prob:code}[{\bf The unrestricted coding problem} ]
     	Let $d,n \in \mathbb N$.
     	Determine 
     	$$
     	\tau_{n,d} :=\inf_{\mathcal C:=\CBRACK{v_i}_{i=1}^n  \subset \mathcal S^{d-1} } \max_{i \neq j} \, \,  \langle v_i, v_j \rangle
     	$$
     \end{prob} 
           \begin{rem}
        	As with preceding problems, the symbol $\tau_{n,d}$ is reserved for solutions to this problem.
        \end{rem}
     
     
        As early as 1955, Rankin~\cite{Rankin1955} provided a perfect solution for the {\it unrestricted coding problem} for all dimensions $d,n \in \mathbb N$ satisfying $1 \leq n \leq 2d$,  providing sharp bounds, replete with examples of optimizers: namely the vertices of simplices and (partial) orthoplexes.
        
        \begin{rem}
        	{\bf Orthoplex} is increasingly used term (particularly within the packing community~\cite{sloane_table, bgb15}) to refer to the vectors corresponding to the vertices of an $\ell^1$-ball;  in other words, an orthoplex is an orthonormal basis unioned with its antipodes.
        \end{rem}

        \begin{thm}[\cite{Rankin1955}, \cite{cwh:1982}; see also \cite{Welch1974}]
        	Let $d,n \in \mathbb N$.
        	\begin{itemize}
        		\item {\bf Simplex Bound}: If $1 \le n \leq d+1$, then $\tau_{n,d} = - \frac{1}{n-1}$ and any $n$-code, $\mathcal C\subset \mathcal S^{d-1}$, that achieves this bound is necessarily a regular ($n-1$)-simplex, meaning all pairwise inner products among the code's elements equal  $\tau_{n,d}$ and $\dim\PARENTH{\SPAN\{\mathcal C\}}=n-1$.
        		\item {\bf Orthoplex Bound}: If $ d+1< n \le 2d$, then $\tau_{n,d} =0$  and any $n$-code, $\mathcal C\subset \mathcal S^{d-1}$, that achieves this bound necessarily contains at least two orthogonal vectors, and all other pairwise inner products occurring within the code are non-positive.  
        		Moreover, if $n=2d$, then the code is necessarily a complete orthoplex.
        		\item{\bf Beyond the orthoplex range} If $n>2d$, then $\tau_{n,d} > 0$ and - except for a few sporadic instances~\cite{lev:2017, 2017arXiv170701858F} -- little is known other than computer-assisted punitively optimal examples, for example, those tabulated in Sloane's online tables of chordal packings~\cite{sloane_table}.
        	\end{itemize}
        \end{thm}
    

     If solutions to the restricted and unrestricted coding problem  equate, ie $\sigma_{n,l, \mathbb F^m} = \tau_{n,\DF}$, then Equation~\ref{eq:gensol} yields the analytic solution to the chordal packing problem for the given parameters via the ``traceless lifting''.   Whenever a given parameter quadruple $(n,l,m, \mathbb F)$, admits a simplex as a solution to both problems, the traceless identity implicates the existence of an {\bf equiangular} tight fusion frame, where a packing is {\bf equiangular} if the set of pairwise trace inner products between its elements is a singleton.    
     
     Of particular note, for the case $l=1$, such objects are more commonly referred to as {\bf equiangular tight frames (ETFs)} and are probably the most famous and well-studied class of optimally spread packings~\cite{2017arXiv170301786H,MR3289408,MR2680063,MR2455575,MR3022574,Boumed_2014,Boumed_2014x,2016arXiv160203490F,Fickus:2015aa,MR2890902,MR2921716,MR3150919,Haas1707:Structures,MR2277977,MR2438918,5946867,Neumaier1989,MR2711357,SustikTroppDhillonHeath2007,Szollosi2013,Szollosi2014b}. Indeed, numerous infinite families are known to exist~\cite{Fickus:2015aa} and dozens~\cite{2017arXiv170301786H,MR3289408,MR2680063,MR2455575,MR3022574,Boumed_2014,Boumed_2014x,2016arXiv160203490F,Fickus:2015aa,MR2890902,MR2921716,MR3150919,Haas1707:Structures,MR2277977,MR2438918,5946867,Neumaier1989,MR2711357,SustikTroppDhillonHeath2007,Szollosi2013,Szollosi2014b} -- if not hundreds -- of mathematicians have contributed to this study, many of whom (see~\cite{Fickus:2015aa}) are engaged in ETF research concurrently with the preparation of this document.   
              
    As for the orthoplex bound, again with $l=1$ fixed, at least three infinite families~\cite{bgb15, WoottersFields1989, BoykinSitharamTarafiWocjan} of tight frames and a few sporadic instances~\cite{Hoggar1982} of parameter quadruples are known to exist where solutions to Problem~\ref{prob:restr}  and Problem~\ref{prob:code} coincide at the orthoplex bound, meaning the conditions on the cardinality range, $n > \DF+1$ and the coherence,
    $\sigma_{n,1, \mathbb F^m} = \tau_{n,\DF}=0$, are satisfied.  Included among these three families are {\it maximal sets of mutually unbiased bases (MUBs)} - to be discussed in further detail shortly, a key ingredient for the constructions of optimally spread packings we are building toward.
    
    Unfortunately, by the inherent restriction of Problem~\ref{prob:restr}, solutions to the two coding problems do not always coincide,  thereby preventing us from ``pulling'' Rankin's solutions back through the traceless identity; otherwise, ETFs would always exist for appropriate parameters.  For example, it is known~\cite{BK06} that an ETF consisting of five elements in $\mathbb R^3$ cannot exist, implicating the inequality  $\sigma_{5,1, \mathbb R^3} > \tau_{5,d_{\mathbb R^3}} = -\frac 1 4$.  Similarly, the nonexistence of an ETF of eight elements in $\mathbb C^3$ was recently verified~\cite{Szollosi2014b}, enforcing the inequality $\sigma_{8,1, \mathbb C^3} > \tau_{8,d_{\mathbb C^3}} = -\frac 1 7$. For more information, we recommend the living table of known ETFs and their properties, located at ~\cite{Fickus:2015aa}.
    
    More generally, for all parameter quadruples, the relationship between solutions to the restricted and unrestricted coding problem may be described succinctly by the inequality, $$\sigma_{n,l, \mathbb F^m} \ge \tau_{n,\DF}.$$  Substituting this inequality into Equation~\ref{eq:gensol} yields the corresponding lower bounds for the chordal coherence of packings for suitable parameter sets.

\begin{thm}\label{th:cohbds}
	[{\cite{Rankin1955}, \cite{cwh:1982}; see also \cite{Welch1974} and \cite{2016arXiv160704546B}}]
	\begin{enumerate}
	\item {\bf Simplex bound}:
	If $\CAL{P}$ is an $(n,l)$-packing for $\mathbb F^m$,  
	then 
	$$
	\mu \PARENTH{\mathcal P} \geq \frac{nl^2-ml}{m(n-1)},
	$$
	and equality is achieved if and only if the fusion frame is { equiangular} and tight.
	\item {\bf Orthoplex bound}:
	If $\CAL{P}$ is an $(n,l)$-packing for $\mathbb F^m$ and  $n >  \DF+1$, then
	$$
	\mu \PARENTH{\mathcal P} \geq \frac{l^2}{m},
	$$
	and if equality is achieved then $\mathcal P$ is optimally spread and  $n \le 2 \DF$.
\end{enumerate}

\end{thm}

We re-iterate that optimally spread fusion frames characterized by the simplex bound are called {\bf equiangular tight fusion frames (ETFFs)}.  As the name suggests, they are indeed both equiangular and tight.  Equiangularity follows immediately from the fact that the $l$-traceless map embeds a given ETFF  into a regular simplex, which is itself equiangular.  Moreover, regular simplexes are zero summing; if $\mathcal P$ is an equiangular, tight $(n,l)$-fusion frame for $\mathbb F^m$ and $\CBRACK{v_{P_i}^{(l)}}_{i=1}^n$ denotes the embedded simplex under $T_l$, then an application of the inverse of the vectorization isomorphism, $\mathcal V^{-1}$, applied to  the overall sum of the embedded, rescaled vectors yields the de-vectorized identity,
\begin{equation}\label{eq:zerotight}\small
{\bf 0_{m\times m}} = \mathcal V^{-1}({\bf 0_{\DF}})= \mathcal V^{-1}\PARENTH{\sum_{i=1}^n v_{P_i}^{(l)} }=
\sqrt{\frac{m}{lm-l^2}} \sum_{i=1}^n\PARENTH{P_i - \frac l m I_m}
=
\sqrt{\frac{m}{lm-l^2}}
\sum_{i=1}^n \PARENTH{P_i - \frac{nl}{m} I_m},\end{equation}
where $\sqrt{\frac{m}{lm-l^2}}$ is the rescaling factor, 
proving the tightness of all ETFFs according to Definition~\ref{def:frames}.

More generally, if $\mathcal P \subset \GRMAN{l}{F}{m}$ is any  $(n,l)$-fusion frame,  $n>\DF+1$ and $\mu\PARENTH{\mathcal P}= \frac{l^2}{m}$, then $\mathcal P$ is optimally spread according to the preceding theorem, in particular characterized by the orthoplex bound.  We call such objects {\bf orthoplectic Grassmannian (ie, optimally spread) $(n,l)$-fusion frames} or {\bf $(n,l)$-OGFFs} for $\mathbb F^m$, or simply {\bf $n$-OGFs} when $l=1$.  In general, OGFFs are not tight; for a thorough examination of this phenomenon, see \cite{2016arXiv160502012C}.  Nevertheless, the three infinite families of OGFs referred to in~\cite{bgb15, WoottersFields1989} are tight and all of the OGFFs to be constructed in this paper are tight.  Finally, we call an OGFF with $n=2\DF$ a {\bf maximal OGFF}, because its tracelessly embedded vectors form a full orthoplex, which is zero-summing, thereby implying tightness according to the same argument for ETFFs; that is,  Equation~\ref{eq:zerotight} applies for maximal OGFFs.

With lower coherence bounds established, we commit the next two subsections to the development of a recipe, previously presented in~\cite{2016arXiv160704546B}, which generates tight OGFFs when certain criteria are satisfied.  
The recipe depends on two major ingredients: mutually unbiased bases and block designs.

\subsection{Two ingredients}
  Because notation becomes somewhat cumbersome as we move forward, when given a natural number, $m$, we will occasionally denote the corresponding index set as follows:
$$
\INDEX{m}:= \CBRACK{1,2,...,m}.
$$

\subsubsection{Mutually Unbiased Bases}
Motivated by consistency and, again, to mitigate notational and typographical issues that will arise in the next section, we make a slight modification to the usual definition~\cite{BandyopadhyayBoykinRoychowdhuryVatan2002, BoykinSitharamTarafiWocjan, GodsilRoy2009, Jaming2010, 1523643, Rao2010, WoottersFields1989} of mutually unbiased bases -- they are usually defined as sets of orthonormal bases with a special property. Equivalently, we reformulate them as families of tight fusion frames comprised of rank one projectors. 
\begin{defn}
	If $\CAL{\Pi}_1=\CBRACK{\pi_j^{(1)}}_{j=1}^m$ and $\CAL{\Pi}_2=\CBRACK{\pi_j^{(2)}}_{j=1}^m$ are a pair of tight $(m,1)$-fusion frames for $\BB{F}^m$, then they are {\bf mutually unbiased} if for $j,k \in \INDEX{m}$ and $s,t \in \INDEX{2}$, the trace inner products satisfy
	$$
	\TR\PARENTH{
		\pi_j^{(s)}
		\pi^{(t)}_k
	} = 
	\left\{\begin{array}{cc} 
	     \frac 1 m,         & s\neq t\\
	     0,                 & s=t, j \neq k\\
	     1,                 & s=t, j=k
	     \end{array}\right..
	$$
  A family of pairwise mutually unbiased bases are simply called {\bf mutually unbiased bases}, or {\bf MUBs}.
\end{defn}

\begin{rem}
	To be clear, by Pareseval's identity, the elements of any tight $(m,1)$-fusion frame, $\CAL{\Pi}=\CBRACK{\pi_j}_{j=1}^m$, for $\mathbb F^m$  (ie, noting $m=n$ here) must arises from an orthonormal basis for $\mathbb F^m$, because the tightness property in Definition~\ref{def:frames} reduces to $$\sum_{i=1}^m \pi_i = I_m.$$
\end{rem}
\begin{notn}
	In this vein, the symbol, $\Pi$, will always refer to a tight $(m,1)$-fusion frame for $\mathbb F^m$ arising from an orthonormal basis.
\end{notn}

Referring back to the orthoplex bound from Theorem~\ref{th:cohbds}, the coherence of a family MUBs  equates with the orthoplex bound when $l=1$, so a sufficiently large family of MUBs can form a Grassmannian frame by the theorem. However, there are many cases where the number of MUBs in $\mathbb F^m$ is known or conjectured~\cite{BoykinSitharamTarafiWocjan,Jaming2010} to be too small to satisfy the cardinality requirement for the orthoplex bound in Theorem~\ref{th:cohbds}. Of seeming relevance is that the number of MUBs is bounded above in terms of the ambient dimension and underlying field.

\begin{thm}\label{th_mub_bd} [Delsarte, Goethals and Seidel \cite{DelsarteGoethalsSeidel1975}] \label{thm:DGS}	
If
	$\{\CAL{\Pi}_j \}_{j=1}^k$ is a family of MUBs for $\BB{F}^m$, then
	$$
	k \le \frac{[\mathbb F : \mathbb R]}{2} m + 1 .
	$$
\end{thm}

Unfortunately, our recipe benefits most greatly from the existence of large families of MUBs, but we encounter frequent deficiencies.
For example, this is especially evident in the real case, because most dimensions, $m$, admit no more than three MUBs in $\mathbb R^m$~\cite{BoykinSitharamTarafiWocjan}; fortunately, whenever $m$ is a power of four~\cite{CameronSeidel1973}, 
real  families of MUBs exist that achieve the upper bound in Theorem~\ref{th_mub_bd}. 
The complex case is a little less severe, as the MUBs' upper cardinality bound is achieved whenever the dimension, $m$, is a prime power \cite{WoottersFields1989}, although -- even for the complex case -- evidence~\cite{Jaming2010} suggests that $7$ MUBs likely do not exist in $\mathbb C^6$, the first complex vector space of composite dimension.

\begin{thm}\label{th_prime_mubs_exist}[\cite{CameronSeidel1973, WoottersFields1989}]\,
	\\
	If $m$ is a prime power, then a family of $m+1$ pairwise mutually unbiased bases for $\BB{C}^m$ exists.
	\\
	If $m$ is a power of 4, then a family of $m/2+1$ pairwise mutually unbiased bases exists for $\BB{R}^m$.
\end{thm}

With Theorem~\ref{th_mub_bd} and Theorem~\ref{th_prime_mubs_exist} in mind, we abbreviate $k_{\mathbb R^m} := m/2+1$ and $k_{\mathbb C^m} := m+1$, and say a family of $k$ MUBs in $\mathbb F^m$ is {\bf maximal} if $k=k_{\mathbb F^m}$.
%

The second ingredient of the recipe are  {\it block designs}.

\subsubsection{Block designs}
As a further effort to avoid convoluted notation and also to illuminate the role of block designs in our construction of tight OGFFs, we diverge somewhat from the conventional symbology found in most combinatorial literature.  Typically, the symbols $v, b, r, k$ and $\lambda$ are designated as the parameters for specifying a given block $t$-design.  To clarify, in the following definition, the symbol $v$ is supplanted by $m$ and $l$ replaces $k$.
\begin{defn}
A {\bf $t$-$(m,l, \lambda)$ block design}, is a pair $\PARENTH{\INDEX{m}, \BLOCKSET}$, where $\BLOCKSET$, is a  collection of subsets of $\INDEX{m}$, called {\bf blocks}, where each block $\CAL{J} \in \BB S$ has cardinality $l$, the cardinality of $\BLOCKSET$, or number of blocks, is $b=\left|\BLOCKSET\right|$, each element of $\INDEX{m}$ occurs in exactly $r$ blocks, and
such that every subset of $\INDEX{m}$ with cardinality $t$ is contained in exactly $\lambda$ blocks.
When the parameters are not important or implied by the context, then $\BLOCKSET$ is also referred to as a {\bf $t$-block design}.
\end{defn}
\begin{notn}
From here on, the pair of symbols $\PARENTH{\INDEX{m}, \BLOCKSET}$ refers to a $t$-block design, where the design's - often suppressed - parameters are prescribed as above.
\end{notn}

A few simple facts about block designs are collected below.
\begin{prop}\label{block:triv}
	Any such block design satisfies the following conditions:
	\begin{enumerate}
		\item $mr=bl$, and
		\item $r(l-1)= \lambda(m-1)$.
	\end{enumerate} 
    Furthermore, it is immediate that for any $t>1$, a $t$-block design is also a $(t-1)$-block design.
\end{prop}

Although we will exploit $t$-block designs with many values for $t$, the forthcoming construction only depends on the existence of certain $1$-block designs, so, henceforth, we will regard all $t$-block designs as $1$-block designs. 
With the basic facts about MUBs and block designs surmised, we are ready to lay out the recipe for tight OGFFs.

\subsection{A recipe for tight orthoplectic Grassmannian fusion frames}

 The key to the recipe is choosing a set of MUBs and a block design whose respective parameters align in such a way that each of the design's blocks inform us on how to choose $l$ elements from a given MUB, which are then summed to form a rank $l$ projector,  constituting one of the $n$ elements of the tight $(n,l)$-OGFF over $\mathbb F^m$ to be constructed; critical to this ``alignment'' of parameters is that the procedure produces sufficiently many projections in order for the orthoplex bound  to apply (ie, we need $n>\DF+1$) and that coherence ultimately equals the orthoplex bound (ie, we also need $\mu(\mathcal P)=\frac{l^2}{m}$).

 In pursuit of the  idea of using a $1$-design's blocks to select rank one projections for us, we define the {\it $\mathcal J$-block projection} as follows.

\begin{defn}
	Given a tight $(m,1)$-fusion frame, $\CAL{\Pi}=\CBRACK{\pi_j}_{j=1}^m$,  for $\BB{F}^m$,  a $1$-$(m,l,\lambda)$ design, and a block  $\CAL J \subset \BLOCKSET$, then the  {\bf $\CAL J$-block projection} with respect to $\CAL{\Pi}$ is
	$$P_{\CAL J} =\sum\limits_{j \in \CAL J} \pi_j.$$
	In this case, writing  $\BLOCKSET=\{\CAL J_1, \CAL J_2, \dots, \CAL J_b\}$, then the $n:=b=\left|\BLOCKSET\right|$-packing,
	 $$\mathcal P := \CBRACK{ P_{\mathcal J_j} }_{j=1}^b \subset \GRASS,$$
	 is called the {\bf block packing} with respect to $\Pi$.
\end{defn}

\begin{prop}\label{prop_1designs}
Given a tight $(m,1)$-fusion frame, $\CAL{\Pi}=\CBRACK{\pi_j}_{j=1}^m$  for $\BB{F}^m$ and a  $1$-$(m,l, \lambda)$ block design where $\BLOCKSET=\{\CAL J_1, \CAL J_2, \dots, \CAL J_b\}$, 
then the block packing with respect to $\Pi$ forms a tight $(n,l,m)$-fusion frame for $\mathbb F^m$, where $n=b$.
\end{prop}
\begin{proof}
	By the preceding definition, every projection formed in this way has rank equal to $l$ and, because $\BLOCKSET$ is a $1$-design, every singleton $\{j\} \subset \INDEX{m}$ occurs exactly $\lambda$ times among the  design's blocks.  Thus, 
	$$
	\sum_{j=1}^b  P_{\mathcal J_j} = 
	\lambda \sum_{j=1}^m  \pi_j = \lambda I_m.
	$$	
\end{proof}

The following fact is perhaps the most crucial and surprising aspect to the construction.
Given a pair of mutually unbiased bases $\mathbb F^m$ and a suitable block design,  one can select  block projections from the respective MUBs to form a pair of tight $(b,l,m)$-fusion frames for $\mathbb F^m$ that achieve the orthoplex bound pair-wise.

\begin{prop}\label{prop_MUBs_give_MUFFs}
	Given $\CAL{\Pi}_i=\CBRACK{\pi_j^{(i)}}_{j=1}^m$, $i \in \INDEX{2}$, a pair of MUBs for $\BB{F}^m$, and $\CAL J, \CAL{J}' \subset \BLOCKSET$ from a $1$-$(m,l,\lambda)$ block design,  then
	$$
	\TR\PARENTH{P_{\CAL J} {P'}_{\CAL J'}} = \frac{l^2}{m},$$
	where
	$P_{\CAL J}$ is the $\CAL J$-block projection with respect to $\CAL{\Pi}_1$ and  ${P'}_{\CAL J'}$ is $\CAL J'$-block projection for ${\CAL{\Pi}_2}.$ 
\end{prop}
\begin{proof}
We compute
$$
\TR(P_{\CAL J} {P'}_{\CAL J'}) 
= \sum\limits_{j \in \CAL J}  \sum\limits_{j' \in \CAL J'} \TR
\PARENTH{
\pi_j^{(1)} \pi_{j'}^{(2)}
}
=  \frac 1 m \sum\limits_{j \in \CAL J}  \sum\limits_{j' \in \CAL J'} 1 = \frac{l^2}{m}
.
$$
\end{proof}

With Propositions~\ref{prop_1designs} and \ref{prop_MUBs_give_MUFFs} in mind,  we see that it is fairly simple to form families of tight fusion frames which achieve the orthoplex  bound pairwise.
However, in order to form optimally spread packings, we must satisfy the sufficiency conditions for the orthoplex bound from Theorem~\ref{th:cohbds}. 
To reiterate, we need {\bf (a)}
$n > \DF +1$ projections in our packing and we need  {\bf (b)} the internal coherence of each block packing arising from each MUB to be sufficiently low.  

To facilitate condition {\bf (b)}, we follow~\cite{2016arXiv160704546B}, recalling their notion of a block design's {\it cohesion}.

\begin{defn}
	Let $(\BB{X}, \BLOCKSET)$ be a $1$-$(m,l,\lambda)$ design. We say that $(\BB{X}, \BLOCKSET)$ is {\bf $c$-cohesive} if there exists $c>0$ such that
	$$
	\max\limits_{\substack{\CAL J, \CAL{J}' \in  \BLOCKSET \\ \CAL J \neq \CAL{J}' }} \left| \CAL J \cap \CAL{J}' \right| \leq c.
	$$
\end{defn}

Finally, the desired recipe for tight OGFFs, originally presented in ~\cite{2016arXiv160704546B}, is described in the following theorem.

\begin{thm}\label{th_costruct_OGFF}
	Let $\left(\INDEX{m},  \BLOCKSET\right)$ be an $l^2/m$-cohesive $1$-$(m,l,\lambda)$ design, where $b:=\left| \BLOCKSET \right|$  and let $\{ \mathcal{\Pi}_k \}_{k \in \INDEX{K}}$ be a set of MUBs for $\BB{F}^m$, where $K b > \DF+1$ and write  $\CAL{\Pi}_k=\CBRACK{\pi_j^{(k)}}_{j=1}^m$ for each $k \in \INDEX{K}$.  If $P_{\CAL J}^{(k)}$ denotes the $\CAL J$-block projection with respect to $\Pi_k$, then the set 
	$$
	\CAL{P} = \bigcup_{k \in  \INDEX{K} } \CBRACK{ P_{\CAL J}^{(k)} }_{\CAL J \in  \BLOCKSET}
	$$
	forms a tight OGFF for $\mathbb F^m$ consisting of $n$ rank $l$ projectors, where
	$
	n =  Kb
	.
	$
\end{thm}

\begin{proof}
By Proposition~\ref{prop_1designs}, $\CBRACK{ P_{\CAL J}^{(k)} }_{\CAL J \in  \BLOCKSET}$ is a tight fusion frame for each $k \in \INDEX{K}$, so $\CAL{P}$ is also a tight fusion frame. 
The cardinality requirement is satisfied since $n > \DF +1$.
Let  $k, k' \in  \INDEX{K}$.
If $k \neq k'$, then 
$$
\TR\PARENTH{ P_{\CAL J}^{(k)} P_{\CAL{J}'}^{(k')}   } = \frac{l^2}{m}
$$
for every $\CAL J, \CAL{J}' \in  \BLOCKSET$ by 
Proposition~\ref{prop_MUBs_give_MUFFs}.
If $k = k'$, then the fact that $( \INDEX{m},  \BLOCKSET)$ is an $l^2/m$-cohesive design yields
\begin{align*}
\max\limits_{\substack{ \CAL J, \CAL{J}' \in  \BLOCKSET \\ \CAL J \neq \CAL{J}'}}\TR\PARENTH{ P_{\CAL{J}}^{(k)} P_{\CAL{J}'}^{(k)}   }
&=
\max\limits_{\substack{ \CAL J, \CAL{J}' \in  \BLOCKSET \\ \CAL J \neq \CAL{J}'}} \sum\limits_{\substack{ j \in \CAL J \\ j' \in \CAL{J}'}}\TR\PARENTH{ \pi_j^{(k)} \pi_{j'}^{(k)} }\\
&=
\max\limits_{\substack{ \CAL J, \CAL{J}' \in  \BLOCKSET \\ \CAL J \neq \CAL{J}'}}
\left| \CAL J \cap \CAL{J}' \right|
\\
&\leq
\frac{l^2}{m},
 \end{align*}
which shows that $\CAL{P}$ is a tight OGFF, as characterized by Theorem~\ref{th:cohbds}.
\end{proof}

Altogether, with respect to the aforementioned {\bf (a)} cardinality issue in mind, it is natural to begin with a maximal set of MUBs, and then seek a compatible  $l^2/m$-cohesive $1$-block design $(\INDEX{m},  \BLOCKSET)$ that produces a sufficient number of coordinate projections per MUB.  This strategy was employed in the original work~\cite{2016arXiv160704546B}, and we similarly depend on the existence maximal MUBs in the following constructions, as well. 
Recall that the {\it known} cases of existence of maximal MUBs are described in Theorem~\ref{th_prime_mubs_exist}.  For each prime power $q$ in the complex case or power of four, $q=4^t$, in the real case, we seek a
$\frac{l^2}{m}$-cohesive $1$-$\PARENTH{m, l, \lambda}$, where $m=q$  and where the number of blocks, $b$, satisfies
\begin{equation}
b > \frac{\DF +1}{\KF}.
\end{equation}

\section{Solutions for the chordal packing problem}\label{sec:sols}

 In~\cite{2016arXiv160704546B}, the authors laid out the preceding construction technique, and then went on to exploit the existence of maximal sets of MUBs in even prime power dimensions and a special class of affine block designs  to construct an infinite family of maximal OGFFs, emphasizing that their construction is highly rigid due to maximality properties of the MUBs and block designs they used, and highlighted the examples'  {cubature properties}. 

In this section, we extend their work, noting that numerous other families of tight OGFFs arise from {\it symmetric block designs} via the final theorem of the preceding section; similarly, we observe examples arising from  other types of affine block designs -- all unmentioned in~\cite{2016arXiv160704546B}.

We begin with symmetric designs, which are  surprising well-suited to this construction principle in a certain sense.
\subsection{Symmetric block designs}
\begin{defn}
	A $t$-$(m,l,\lambda)$-block design is {\bf symmetric} if $m=b$ or, equivalently, if $l=r$. 
\end{defn}

As verified in~\cite{cameron_lint_1991}, symmetric block designs have the useful property that the pairwise block intersections is constant.
\begin{thm}\label{th:sym.eq}[\cite{cameron_lint_1991}]
	For a symmetric $t$-$(m,l,\lambda)$ block design, $\left(\INDEX{m},  \BLOCKSET\right)$,  every $\mathcal J_j, \mathcal J_{j'} \in \BLOCKSET$ with $j \neq j'$ satisfies $$\left| \mathcal J_j \cap \mathcal J_{j'} \right| = \lambda.$$
\end{thm}

\subsubsection{Simple ETFFs from all t-block designs}
An immediate -- albeit trivial -- corollary to this observation is that every symmetric block designs yields an ETFF.

\begin{cor}\label{prop:symmequi}
	Given any symmetric $t$-$(m,l,\lambda)$ block design and a tight  $(m, 1)$-fusion frame, $\Pi = \CBRACK{\pi_i}_{i=1}^m$, for $\mathbb F^m$,  then the family of $n=b$ block projections for $\Pi$ forms an ETFF  for $\mathbb F^m$. 
\end{cor} 
\begin{proof}
	By design, Proposition~\ref{prop_1designs} ensures that the sequence of $n=b$ rank $l$ block projections, $\mathcal P = \CBRACK{P_i := \sum_{j \in \mathcal J_i}  \pi_j}_{i=1}^m$ forms a tight fusion frame for $\mathbb F^m$.  By Theorem~\ref{th:cohbds}, these are chordally equiangular, since, for $i \neq i'$,  we have
	$$
	\TR\PARENTH{P_i P_{i'} } = \sum_{j \in \mathcal J_i}\sum_{j' \in \mathcal J_{i'}} \TR\PARENTH{  \pi_j   \pi_{j'} } = \lambda,
	$$ 
	where one can confirm that $\lambda$ is the packing constant for this special case.
\end{proof}

As implicated by the numerous examples of symmetric block designs in ~\cite{MR2246267} and the references therein,  these simple examples of ETFFs exist in abundance.

\subsubsection{Tight OGFFs from certain block dseigns}

Of more significance is the multitude of examples of tight OGFFs that arise from the construction outlined in theorem~\ref{th_costruct_OGFF}, depending heavily on the existence of maximal MUBs. The surprise here is that every nontrivial symmetric $t$-block-design satisfies the $l^2/m$-cohesive property.
\begin{prop}
	Assume $1 < l < m-1$.  Every symmetric $t$-$(m,l,\lambda)$ block design, $(\INDEX{m},  \BLOCKSET)$, is $l^2/m$-cohesive.
\end{prop}
\begin{proof}
	By Theorem~\ref{th:sym.eq}, we only need to show that $\lambda \le l^2 /m$, which follows by applying the trivial necessary conditions from Proposition~\ref{block:triv} to obtain the elementary estimation,
	$$
	\lambda = \frac{r(l-1)}{m-1} = \frac{l(l-1)}{m-1} \le \frac{l^2}{m}.
	$$
\end{proof}

Thus, following the strategy, wherein we employ maximal MUBs, as outlined after the statement of Theorem~\ref{th_costruct_OGFF}, we seek symmetric block designs for which $m=b=q$, where $q$ is a prime power or a power of four, for the complex and real cases, respectively.

\begin{cor}\label{cor:symogffs}
	Assume $1 < l < m-1$.
	\begin{itemize}
		\item 	Every symmetric $t$-$(m,l,\lambda)$ block design for which $m$ is a prime power yields a tight $(n,l)$-OGFF for $\mathbb C^m$, where $n=m(m+1)$, according to the construction described in Theorem~\ref{th_costruct_OGFF}.
		\item 	Every symmetric $t$-$(m,l,\lambda)$ block design for which $m$ is a power of four yields a tight $(n,l)$-OGFF for $\mathbb R^m$, where $n=m(m/2 +1)$, according to the construction described in Theorem~\ref{th_costruct_OGFF}.
	\end{itemize}
\end{cor}

Examining the list of known families of symmetric designs in II.6.8 of \cite{MR2246267}, we identify numerous families of symmetric block designs for which Corollary~\ref{cor:symogffs} yields more tight OGFFs.

\subsubsection{Point-hyperplane symmetric designs; see Family 1 of~ \cite{MR2246267}}
Given a prime power $q$ and $t \ge 2$, a symmetric $1$-$\PARENTH{\frac{q^{t+1} -1}{q-1}, \frac{q^t-1}{q-1}, \frac{q^{t-1}-1}{q-1} }$ block design exists, known as a {\bf point-hyperplane design}. 
\begin{itemize}
	\item 
	By Corollary~\ref{cor:symogffs}, a tight $\PARENTH{m, \frac{q^t-1}{q-1}}$-OGFF for $\mathbb C^m$, comprised of $n:=\frac{(q^{t+1} -1)(q^{t+1} +q -2)}{(q-1)^2}$ projections,  exists whenever $m:=\frac{q^{t+1} -1}{q-1}$ is a prime power, leading to numerous families of examples of tight OGFFs. 
	\item
	Unfortunately, the necessary conditions of these block designs seem too restrictive to admit real examples of OGFFs according to our approach.
\end{itemize}
Contingent upon the open question regarding the existence of an infinitude of Fermat or Mersenne primes -- a famous open problem in number theory~\cite{MR0327647} -- complex examples arising from this family may be infinite by taking $q=2$.

\subsubsection{Hadamard symmetric designs; see Family 2 of~ \cite{MR2246267}}
Given $t \ge 1$, a symmetric $1$-$\PARENTH{4t-1, 2t-1, t-1 }$ block design exists, known as a {\bf Hadamard design}. 
\begin{itemize}
	\item 
	By Corollary~\ref{cor:symogffs}, a tight $\PARENTH{m, 2t-1}$-OGFF for $\mathbb C^m$, comprised of $n:=4t(4t-1)$ projections,  exists whenever $m:= 4t-1$ is a prime power, leading to an infinite family of complex examples of tight OGFFs.
	\item
	Unfortunately,  $m:= 4t-1$ is never a power of four, so these designs yield no real examples.
\end{itemize}

\subsubsection{Menon symmetric designs; see Family 6 of~ \cite{MR2246267}}
The existence of a Hadamard matrix (see \cite{MR2246267} for details) of order $4t^2$, $t \in\mathbb N$, is equivalent to the existence of a symmetric $t$-$\PARENTH{4t^2, 2t^2-t, t^2-t}$ block design, called a {\bf Menon design}.  According to the well-known Hadamard conjecture~\cite{MR2246267}, such block designs exist for all values of $t$.  Besides the conjecture, a simple tensor construction~\cite{MR2246267}, among other constructions, assures their existence when $t=2^s$ for some $s \in \mathbb N$.
\begin{itemize}
	\item 
	By Corollary~\ref{cor:symogffs}, when $t=2^s$ for some $s \in \mathbb N$, a tight $\PARENTH{4t^2, 2t^2-t}$-OGFF for $\mathbb C^m$ comprised of $n:=16t^4 + 4t^2$ projections exists,  leading to an infinite family of complex examples. 
	\item
	Conveniently, these designs are well-suited for the real case.  Similarly, when $t=2^s$ for some $s \in \mathbb N$, where $s$ is even, then a tight $\PARENTH{4t^2, 2t^2-t}$-OGFF for $\mathbb R^m$ comprised of $n:=4t^2(2t^2 + 1)$ projections exists,  leading to an infinite family of tight OGFFs.
\end{itemize}

\subsubsection{Wallis symmetric designs; see Family 7 of~\cite{MR2246267}}
Given a prime power $q$ and $t \in \mathbb N$, a symmetric 
$t$-$\PARENTH{m, l,\lambda  }$ block design exists, known as a {\bf Wallis design},
where
\begin{itemize}
	\item $m:=q^{t+1}\left(q^t + q^{t-1} + \dots + q^2 + q + 2 \right),$\\
	\item $l:= q^{t}(q^t + q^{t-1} + \dots + q^2 +  q + 1 )$, and \\
	\item $\lambda =  q^{t}(q^{t-1} + q^{t-2} + \dots + q^2 + q + 1 ).$
\end{itemize}	
	  Taking $q=2$
	  and computing (inductively or directly) that
	  $$
	  2^{t+1}(2^t + 2^{t-1} + .... + 4 + 2 + 2) = 4^{t+1},
	  $$
 the Wallis family yields infinite families of both real and complex tight OGFFs.
\begin{itemize}
	\item 
	By Corollary~\ref{cor:symogffs}, when $q=2$ and $t \in \mathbb N$, a tight $\PARENTH
	{n,l}$-OGFF for $\mathbb C^m$ comprised of $n:=4^{2t+2} + 4^{t+1}$ projections exists,  leading to an infinite family of complex examples.
	\item
	Once again, these designs are well-suited for the real case.  When $q=2$ and $t \in \mathbb N$, a tight $\PARENTH{n, l}$-OGFF for $\mathbb R^m$ comprised of $n:=4^{t+1}(\frac{4^{t+1}}{2}+1)$ projections exists,  leading to another infinite family of real examples.
\end{itemize}

\subsubsection{Wilson/Brouwer symmetric designs; see Family 1 of~11 \cite{MR2246267}}
Given an odd prime power $q$ and $t \in \mathbb N$, a symmetric 
$1$-$\PARENTH{m, l,\lambda  }$ block design exists, which we refer to as a {\bf Wilson/Brouwer  design},
where
\begin{itemize}
	\item $m:=2\left(q^t + q^{t-1} + \dots + q^2 + q \right) +1,$\\
	\item $l:= q^{t},$ and \\
	\item $\lambda =  q^{t-1}\frac{q-1}{2}$.
\end{itemize}
Unfortunately, the addition of $1$ to an even number makes it impossible for $m$ to be a power of four here.  Still, for various parameters -- eg, $q=3$, $t=1$ yields $m=7$ -- this family produces various complex instances of tight OGFFs.

A comprehensive list of all tight OGFFs obtained from symmetric designs is beyond the scope of this work -- particularly due to outstanding open questions about their existence for various parameters~\cite{MR2246267}.  Having demonstrated families arising in this fashion, we move to tight OGFFs arising from affine designs.

\subsection{Affine Designs}
In the work where this recipe is presented, the authors combined real and complex even powers of maximal MUBs of Hadamard $3$-designs (see\cite{2016arXiv160704546B}), a special case of affine designs, to form maximal OGFFs according to Theorem~\ref{th_costruct_OGFF}. 

\begin{defn}
A block design $\PARENTH{\INDEX{m}, \BLOCKSET}$  is {\bf resolvable} if  $\BLOCKSET$ partitions into subsets, called {\bf parallel classes}, such that
\begin{itemize}
\item the blocks with in each class are disjoint, and
\item for each parallel class, every element of $\INDEX{m}$ is contained in a block.
\end{itemize} 
Moreover, if
the number of elements occurring in the intersection between blocks from different parallel classes is constant, then it is an {\bf affine design}.
	
\end{defn}

A simple result in~\cite{MR2246267} assures the existence of analogous block designs for the odd prime power case.
We restate these for our special case.

\begin{prop}\label{prop_affexist}[See II.7.10 \cite{MR2246267}]
	If $m=p^{t+1}$, $l=p^t$ and $\lambda=\frac{p^t-1}{p-1}$ for some prime, $p$, then a resolvable $1$-$(m,l,\lambda)$ block design exists and its remaining parameters are therefore  $r=\sum_{i=0}^t p^i$ and  $b=\sum_{i=1}^{t+1} p^i$.
\end{prop}

 According to Bose's condition, such designs are $l^2/m$-cohesive affine block designs.
 
 \begin{thm}[Bose's condition; see Theorem II.7.28 of \cite{MR2246267}]
 	Given any resolvable $1$-$(m,l,\lambda)$ block design, the number of blocks is bounded by the other parameters according to
 	  $$b \ge m +r-1,$$
 	  and this lower bound is achieved if and only if the design is an $l^2/m$-cohesive affine design.
 \end{thm}
 
 One may verify that the designs produced by Proposition~\ref{prop_affexist} achieve this lower bound.  Noting that such designs satisfy
 $$k_{\mathbb C^{m}} b= (p^{t+1} +1) \sum_{i=1}^{t+1} p^i > d_{\mathbb C^m}+1 = p^{2t+2},$$
 the conditions of Theorem~\ref{th_costruct_OGFF} are satisfied.
 
 \begin{cor}
 	For every prime $p$ and $t\in \mathbb N$, a tight $(n,l)$-OGFF, comprised of $n=(p^{t+1} +1) \sum_{i=1}^{t+1} p^i$ exists for $\mathbb C^m$, where $m=p^{t+1}$ and $l=p^t$.
 \end{cor} 

Although these tight OGFFs are not maximal OGFFs, they are maximal with respect to the construction in Theorem~\ref{th_costruct_OGFF}, as each instance is produced by a maximal set of MUBs and a maximal affine design, according to Bose's condition.


\bibliography{universal}
\bibliographystyle{plain}

\end{document}